\documentclass{amsproc}
\usepackage{graphicx}
\usepackage{amscd}
\usepackage{amsmath}
\usepackage{amsfonts}
\usepackage{amssymb}
\newtheorem{theorem}{Theorem}
\theoremstyle{plain}

\newtheorem{corollary}{Corollary}

\newtheorem{definition}{Definition}

\newtheorem{proposition}{Proposition}
\newtheorem{remark}{Remark}

\newtheorem{summary}{Summary}
\numberwithin{equation}{section}

\begin{document}
\title[On Algebraic Properties]{On Algebraic Properties of Sets of Banach Ideal Function Spaces}
\author{Eugene Tokarev}
\address{B.E. Ukrecolan, 33-81, Iskrinskaya str., 61005, Kharkiv-5, Ukraine}
\email{tokarev@univer.kharkov.ua}
\thanks{This paper is in final form and no version of it will be submitted for
publication elsewhere.}
\subjclass{Primary 46E30, 46B70; Secondary 46B10, 46B20}
\keywords{Lattices, Galois connexions, Ideal Banach function spaces}
\dedicatory{Dedicated to the memory of S. Banach.}
\begin{abstract}It is shown that a set $\mathcal{J}\left(  \mu\right)  $ of Banach lattices of
real-valued measurable functions, defined on a measure space $\left(
\Omega,\Sigma,\mu\right)  $, may be equipped with a some natural ordering
under which it becomes a distributive lattice, which is Dedekind complete
provided $\mu$ is a probability measure. Moreover, some natural operations on
considered spaces are in Galois connexion. These results are of most interest
for symmetric Banach spaces.
\end{abstract}
\maketitle

\section{Introduction}

This paper is devoted to study algebraic properties of a set $\mathcal{J}%
\left(  \mu\right)  $ of Banach ideal spaces of real valued $\mu$-measurable
functions. Namely, it will be shown that a quite natural ordering
''$\subset^{1}$'' on $\mathcal{J}\left(  \mu\right)  $ makes this set to be a
lattice; some restrictions on spaces from $\mathcal{J}\left(  \mu\right)  $
mark out sublattices of $\left\langle \mathcal{J}\left(  \mu\right)
,\subset^{1}\right\rangle $ having nice algebraic properties. Compositions of
some natural operations on Banach ideal spaces (such as either the operation
of conversion of a given space $E$ into its dual $E^{\prime}$ or the operation
to pick out all elements of $E$ having an absolutely continuous norm to
generate a new space $E_{0}$) may be chosen in a such way that they will be in
the Galois connexion.

Section 2 is devoted to recall some definitions and notations that touch on
Banach ideal spaces. The commonly used terminology is widely changed from one
paper to another. Below mainly will be used the terminology of reviews [1] and
[2]. For all results that are mentioned below without proofs the reader refers
to these reviews.

Original results are contained in sections 3, 4 and 5..

\section{Definitions}

Let $\left(  \Omega,\Sigma,\mu\right)  $ be a measure space, i.e., an abstract
set $\Omega$ with a $\sigma$-algebra $\Sigma$ of its subsets and a countably
additive function (measure) $\mu$, defined on $\Sigma$ with the range in
$\mathbb{R}_{+}$.

Let $L_{0}\left(  \mu\right)  =L_{0}\left(  \Omega,\Sigma,\mu\right)  $ be the
set of (classes of) $\mu$-measurable real valued functions, defined on
$\Omega$. Certainly, $L_{0}\left(  \mu\right)  $ is a vector space under usual
operations of addition of functions and multiplication by a scalar.
$L_{0}\left(  \mu\right)  $ is also a lattice under a natural partial order
($x\left(  \omega\right)  \leq y\left(  \omega\right)  $ means that $x\left(
\omega\right)  \leq y\left(  \omega\right)  $ a.e.).

An ideal Banach function space (shortly: BIS) $E\left(  \mu\right)  $ is a
vector subspace of $L_{0}\left(  \mu\right)  $, which is equipped with a
\textit{Banach norm} $\left\|  \cdot\right\|  _{E}$ (i.e. $E\left(
\mu\right)  $ is complete in the norm topology), which is\textit{\ monotone},
i.e., such that

\begin{center}
$y\left(  \omega\right)  \in E\left(  \mu\right)  ,$ $x\left(  \omega\right)
\in L_{0}\left(  \mu\right)  $\ \textit{and }$\left|  x\left(  \omega\right)
\right|  \leq\left|  y\left(  \omega\right)  \right|  $\textit{\ implies that}
$\left\|  x\right\|  _{E}\leq\left\|  y\right\|  _{E}$.
\end{center}

Classical examples of BIS are Lebesgue-Riesz spaces $L_{p}(\mu)$ ($1\leq
p\leq\infty$). However there are examples of such measure space $\left(
\Omega,\Sigma,\mu\right)  $ (which does not have the direct sum property;
definitions see below) that the space $L_{p}(\mu)$ is not complete (i.e., is
not a Banach space). Certainly, it may be completed (by the usual procedure of
completition) and the resulting Banach space will be of kind $L_{p}%
(\mu^{\prime})$ as well. However the measure space $\left(  \Omega^{\prime
},\Sigma^{\prime},\mu^{\prime}\right)  $, where the complete space
$\overline{L_{p}(\mu)}=L_{p}(\mu^{\prime})$ will be defined differs from the
initial one.

So, it is necessary to put some restrictions on the measure space $\left(
\Omega,\Sigma,\mu\right)  $.

\begin{definition}
A measure space $\left(  \Omega,\Sigma,\mu\right)  $ is said to be admissible
if is satisfies the following conditions:

\begin{itemize}
\item  If $A\in\Sigma$; $\mu\left(  A\right)  =0$ and $B\subset A$ then
$B\in\Sigma$ and $\mu\left(  B\right)  =0$ (the measure $\mu$ is complete).

\item  If $A\subset\Omega$ and every $B\in\Sigma$ with $\mu(B)<\infty$ is so
that\ $A\cap B\in\Sigma$, then $A\in\Sigma$.

\item  If $A\in\Sigma$ and $\mu\left(  A\right)  =\infty$ then there exists
$B\subset A$, $B\in\Sigma$ such that $\mu(B)<\infty$ (the measure $\mu$ is semifinite).

\item  There exists a set $\{A_{i}:i\in I\}$ of pairwice disjoint subsets of
$\Omega$ with $\mu\left(  A_{i}\right)  <\infty$ for every $i\in I$ so that

\begin{itemize}
\item  Every $B\in\Sigma$ of finite measure $\mu\left(  B\right)  <\infty$ may
be represented as
\[
B=\cup\{B\cap A_{i}:i\in I_{0}\}\cup A_{0},
\]
where $I_{0}$ is a countable subset of $I$ and $\mu\left(  A_{0}\right)  =0$
(the measure $\mu$ has the direct sum property).
\end{itemize}
\end{itemize}
\end{definition}

It is known (see e.g. [1] and [2]) that for every admissible measure space
$\left(  \Omega,\Sigma,\mu\right)  $ each BIS $E\left(  \mu\right)  $ is
conditionally Dedekind complete, and a set of all integral functionals on
$E\left(  \mu\right)  $ is total over it.

For any $A\in\Sigma$ the triple $\left(  A,\Sigma_{A},\mu_{A}\right)  $, where
$\Sigma_{A}=\{B\cap A:B\in\Sigma\}$ and $\mu_{A}\left(  B\right)  =\mu\left(
B\cap A\right)  $ is a restriction of $\mu$ is an admissible set provided
$\left(  \Omega,\Sigma,\mu\right)  $ is admissible.

It will be said that BIS $E\left(  \mu\right)  $ is of \textit{maximal width
in} $\left(  \Omega,\Sigma,\mu\right)  $ if
\[
\{z\in S\left(  \mu\right)  :zy=0\text{ \ for all }y\in E\left(  \mu\right)
\}=0.
\]

\begin{definition}
Let $\left(  \Omega,\Sigma,\mu\right)  $ be an admissible measure space. A set
$\mathcal{J}\left(  \mu\right)  $ is the set of all BIS $E\left(  \mu\right)
$ that are of maximal width in $\left(  \Omega,\Sigma,\mu\right)  $.
\end{definition}

So, $E\left(  \mu_{A}\right)  \in\mathcal{J}\left(  \mu\right)  $ if and only
if $\mu\left(  \Omega\backslash A\right)  =0$.

Let $I=\left\langle I,\ll\right\rangle $ be a partially ordered set. It will
be said that $I$ is a \textit{directed set} if for any $i_{1}$, $i_{2}\in I$
there exists $i\in I$ such that $i_{1}\ll i$ and $i_{2}\ll i$.

A sequence $\{x_{i}:i\in I\}$, which is indexed by elements of the directed
set $I=\left\langle I,\ll\right\rangle $ will be called a\textit{\ net}. It
will be written $x_{i}\downarrow$ if $i\ll j$ implies $x_{i}\geq x_{j}$. If
$x_{i}\downarrow$ and $\inf_{i\in I}\left(  x_{i}\right)  =x_{0}$, we shall
write $x_{i}\downarrow x_{0}$.

Let $E=E\left(  \mu\right)  $ be a BIS; $x\in E$. It will be said that the
norm of $x$ is \textit{order continuous }(shortly: $\left(  o\right)
$-\textit{continuous}) provided the condition $\left|  x\right|  \geq
x_{i}\downarrow0$ implies that $\lim_{I}\left\|  x_{i}\right\|  _{E}=0$.

The set of all elements of $E$ having the $\left(  o\right)  $-continuous norm
is denoted by $E_{0}$.

Certainly, $E_{0}$ is a closed Banach subspace of $E$ and is an \textit{ideal}
in $E$:

\begin{center}
\textit{if }$x\in E_{0}$\textit{; }$y\in E$\textit{\ and }$\left|  y\right|
\leq\left|  x\right|  $\textit{\ then }$y\in E_{0}$.
\end{center}

Notice that $E_{0}$ need not to be of maximal width in $\left(  \Omega
,\Sigma,\mu\right)  $; moreover it may be trivial. E.g., $\left(  L_{\infty
}\left[  0,1\right]  \right)  _{0}=\{0\}$.

Recall that a subset $F$ of a BIS $E$ is said to be a \textit{foundation in
}$E$ if it is an ideal in $E$ and is of maximal width in $\left(
\Omega,\Sigma,\mu\right)  $).

Let $E\left(  \mu\right)  \in\mathcal{J}\left(  \mu\right)  $. A \textit{dual}
space $E^{\prime}$ is the space of all elements $f(t)\in L_{0}(\mu)$ such
that
\[
\left\|  f\right\|  _{E^{\prime}}=\sup\{\int\nolimits_{\Omega}f\left(
t\right)  x\left(  t\right)  d\mu:\left\|  x\right\|  _{E}=1\}<\infty.
\]

$E^{\prime}$ may be identified with a subset of the conjugate space $E^{\ast}
$: every element $f\in E^{\prime}$ generates the (\textit{integral})
functional $f\in E^{\ast}$ by the rule:
\[
\left\langle f,x\right\rangle =\int\nolimits_{\Omega}f\left(  t\right)
x\left(  t\right)  d\mu.
\]

Certainly, $E^{\prime}$ is a Banach space under the norm $\left\|
\cdot\right\|  _{E^{\prime}}$ and is a BIS (of maximal width), which belongs
to $\mathcal{J}\left(  \mu\right)  $.

\begin{remark}
There may be situations when $E_{0}$ is a foundation in $E$ but $\left(
E^{\prime}\right)  _{0}$ is not a foundation in\ $E^{\prime}$. E.g.,
$E=E_{0}=L_{1}\left[  0,1\right]  $; $E^{\prime}=L_{\infty}\left[  0,1\right]
$ and $\left(  E^{\prime}\right)  _{0}=\{0\}$. The paper $[3]$ contains an
example of such BIS $E$ that $E_{0}=\left(  E^{\prime}\right)  _{0}=\{0\}$.
\end{remark}

\begin{definition}
Let $\left(  \Omega,\Sigma,\mu\right)  $ be an admissible measure space. A set
$\mathcal{J}_{0}\left(  \mu\right)  $ is the set of all Banach ideal spaces
$E=E\left(  \mu\right)  $ such that $E_{0}=E_{0}\left(  \mu\right)  $ is a
foundation in $E\left(  \mu\right)  $ and $\left(  E^{\prime}\right)  _{0} $
is a foundation in $E^{\prime}$.
\end{definition}

Let $\left(  \Omega,\Sigma,\mu\right)  $ be an admissible measure space.

Let $E\subset F$ be BIS. Define an operator $i\left(  E,F\right)
:E\rightarrow F$, which asserts to any $x\in E$ the same function $x\in F$.
This operator is called \textit{the inclusion operator}. Its norm (that is the
infimum of all possible constants $c(E,F)$) is called \textit{the inclusion
constant}.

The relation $E\subset^{1}F$\ means that the inclusion constant $c(E,F)=1$.

The class $\mathcal{J}\left(  \mu\right)  $ is partially ordered by the
relation $E\subset^{1}F$

\begin{definition}
$\left\langle \mathcal{J}\left(  \mu\right)  ,\subset^{1}\right\rangle $ is a
partially ordered set, in which two BIS $E$ and $F$ are identified if and only
if $E\subset^{1}F$ and $F\subset^{1}E$. This means that $E$ and $F$ are
identical as sets, as vector lattices, as topological vector spaces and as
Banach spaces.
\end{definition}

On $\mathcal{J}\left(  \mu\right)  $ the more general relation $E\subset^{c}F
$ may be defined. It means that the inclusion constant $c(E,F)$ is bounded.
This relation partially orders the quotient set $\mathcal{J}^{\approx}\left(
\mu\right)  =\mathcal{J}\left(  \mu\right)  /\approx$, where the equivalence
relation $E\approx F$ means that $E\subset^{c}F$ and $F\subset^{c}E$.

Below it will be regarded only the partially ordered set $\left\langle
\mathcal{J}\left(  \mu\right)  ,\subset^{1}\right\rangle $. The study of the
set $\left\langle \mathcal{J}^{\approx}\left(  \mu\right)  ,\subset
^{c}\right\rangle $ is reserved to readers.

\begin{remark}
It is worthwhile to note that the definition of relations $\subset^{1}$ and
$\subset^{c}$ in a general case has some defects (which are eliminated while
symmetric Banach spaces will be considered).

E.g., spaces $E=L_{p}\left[  0,1/2\right]  \oplus_{1}L_{q}\left[
1/2,1\right]  $ and $F=L_{q}\left[  0,1/2\right]  \oplus_{1}L_{p}\left[
1/2,1\right]  $ \ for $p\neq q$ must be considered as different. Indeed, they
are not compatible neither in the sense of $\subset^{1}$ nor in the sense of
$\subset^{c}$.
\end{remark}

The following summary assembles all results about the order $\subset^{1}$ and
operations $E\rightsquigarrow E_{0}$; $E\rightsquigarrow E^{\prime}$ that will
be needed later.

\begin{summary}
Let $\left(  \Omega,\Sigma,\mu\right)  $ be an admissible set; $\mathcal{J}%
\left(  \mu\right)  $ - the corresponding set of Banach ideal spaces of $\mu
$-measurable functions defined on $\Omega$ that are of maximal width. Let $E$,
$F\in\mathcal{J}\left(  \mu\right)  $. Then

\begin{enumerate}
\item $E_{0}\subset^{1}E$;

\item $E\subset^{1}F$ implies that $E_{0}\subset^{1}F_{0}$;

\item $\left(  E_{0}\right)  _{0}=E_{0}$;

\item $E\subset^{1}E^{\prime\prime}$;

\item $E^{\prime}=E^{\prime\prime\prime}$;

\item  If $E\subset^{1}F$ then $F^{\prime}\subset^{1}E^{\prime}$.
\end{enumerate}
\end{summary}

\section{$\left\langle \mathcal{J}\left(  \mu\right)  ,\subset^{1}%
\right\rangle $ as a lattice}

Consider an admissible measure space $\left(  \Omega,\Sigma,\mu\right)  $ and
the corresponding set $\mathcal{J}\left(  \mu\right)  $ equipped with the
partial order $\subset^{1}$.

Let $E$, $F\in\mathcal{J}\left(  \mu\right)  $. According to [4] define a pair
of BIS: $E\cap F$ and $E+F$.

$E\cap F$ consists of all functions $x$, that are common to $E$ and $F$: $f\in
E\cap F$ and is equipped with the norm
\[
\left\|  f\right\|  _{E\cap F}=\max\{\left\|  f\right\|  _{E},\left\|
f\right\|  _{F}\}.
\]

$E+F$ is formed by functions of kind $f=u+v$, $u\in E$; $v\in F$, such that
\[
\left\|  f\right\|  _{E+F}=\inf\{\left\|  u\right\|  _{E}+\left\|  v\right\|
_{F}:u+v=f\}<\infty;
\]

Define on $\left\langle \mathcal{J}\left(  \mu\right)  ,\subset^{1}%
\right\rangle $ binary operations, $\vee$ and $\wedge$. Namely, put
\[
E\vee F:=E+F;\text{ \ \ }E\wedge F:=E\cap F,
\]

\begin{theorem}
$\left\langle \mathcal{J}\left(  \mu\right)  ,\vee,\wedge\right\rangle $ is a lattice.
\end{theorem}

\begin{proof}
Recall that an algebraic structure $A$, endowed with a pair $\vee,\wedge$ of
binary mappings
\[
\vee:A^{2}\rightarrow A;\text{ \ }\wedge:A^{2}\rightarrow A
\]
is a lattice, if for any $a$, $b$, $c\in A$
\begin{align*}
a\vee b &  =b\vee a;\text{ \ }a\wedge b=b\wedge a\\
a\vee\left(  b\vee c\right)   &  =\left(  a\vee b\right)  \vee c;\text{
\ }a\wedge\left(  b\wedge c\right)  =\left(  a\wedge b\right)  \wedge c\\
a\vee\left(  a\wedge b\right)   &  =a\wedge\left(  a\vee b\right)  =a.
\end{align*}
It is an easy exercise to verify lattice axioms for $\left\langle
\mathcal{J}\left(  \mu\right)  ,\vee,\wedge\right\rangle $.
\end{proof}

\begin{remark}
In a general case $\left\langle \mathcal{J}\left(  \mu\right)  ,\subset
^{1}\right\rangle $ is not Dedekind complete.

Indeed, consider a BIS $E$ and a sequence $\left\langle E,\left\|
\cdot\right\|  _{n}\right\rangle _{n\in\mathbb{N}}$, where $\left\|
x\right\|  _{n}=n\left\|  x\right\|  $ for all $n\in\mathbb{N}$. Certainly,
such a sequence does not have any upper bound.
\end{remark}

This difficulty may be overcame if $\mu$ is a probability measure, say,
$\mathbb{P}$, which does not contain atoms of positive measure.

Let $\mathcal{J}\left(  \mathbb{P}\right)  $ be a set of all BIS that
satisfies the norming condition $\left\|  \chi_{\Omega}\left(  t\right)
\right\|  =1$, where $\chi_{A}\left(  t\right)  $ is the indicator function of
$A\in\Sigma$:$\ $%
\[
\chi_{A}\left(  t\right)  =1\text{ \ for \ }t\in A;\text{ \ }\chi_{A}\left(
t\right)  =0\text{ \ for \ }t\notin A
\]
and are of maximal width on the probability space. The following theorem is valid.

\begin{theorem}
The lattice $\left\langle \mathcal{J}\left(  \mathbb{P}\right)  ,\vee
,\wedge\right\rangle $ is Dedekind complete.
\end{theorem}

\begin{proof}
Indeed, if $E\left(  \mathbb{P}\right)  $ satisfies the norming conditions
then
\[
L_{\infty}\left(  \mathbb{P}\right)  \subset^{1}E\left(  \mathbb{P}\right)
\subset^{1}L_{1}\left(  \mathbb{P}\right)  .
\]
The greatest lower bound of a family $\left(  E_{i}\right)  _{i\in I}$ of BIS
that satisfy norming conditions is the space $\Cap E_{i}$, which consists of
all elements that are common to all $E_{i}$'s. Its norm is given by
\[
\left\|  x\right\|  _{\Cap E_{i}}=\sup\{\left\|  x\right\|  _{E_{i}}:i\in I\}.
\]
Since every set of BIS that satisfy the norming condition is bounded (by
$L_{1}\left(  \mathbb{P}\right)  $), the least upper bound of $\left(
E_{i}\right)  _{i\in I}$ may be obtained as the greatest lower bound of a
family of all upper bounds of\ $\left(  E_{i}\right)  _{i\in I}$.

Another way is to define (according to $[4]$) a space $\Cup E_{i}$ as an ideal
in $L_{0}\left(  \mathbb{P}\right)  $ that consists of all $x\in L_{0}\left(
\mathbb{P}\right)  $, which has the representation of kind
\[
x=\sum\nolimits_{i\in I}u_{i}\text{ \ where \ }u_{i}\in E_{i};\text{ \ }%
\sum\nolimits_{i\in I}\left\|  u_{i}\right\|  _{E_{i}}<\infty.
\]
It is equipped with the norm
\[
\left\|  x\right\|  _{\Cup E_{i}}=\inf\{\sum\nolimits_{i\in I}\left\|
u_{i}\right\|  _{E_{i}}:x=\sum\nolimits_{i\in I}u_{i};\text{ \ }u_{i}\in
E_{i}\}.
\]
From results of $[4]$ it follows that $\sup\nolimits_{i\in I}\left(
E_{i}\right)  =\Cup E_{i}$.
\end{proof}

Now return to an arbitrary admissible measure space $\left(  \Omega,\Sigma
,\mu\right)  $ and consider the lattice $\left\langle \mathcal{J}\left(
\mu\right)  ,\vee,\wedge\right\rangle $.

Let $E$, $F\in\mathcal{J}\left(  \mu\right)  $; $E\subset^{1}F$. An order
interval $[E,F]_{\subset^{1}}$ is given by
\[
\lbrack E,F]_{\subset^{1}}:=\{H\in\mathcal{J}\left(  \mu\right)  :E\subset
^{1}H\subset^{1}F\}.
\]
It will be denoted by $[E,F]$ for simplicity. For any order interval $[E,F]$
there may be defined mappings $\lambda=\lambda_{E,F}$ and $\rho=\rho_{E,F}$ as
follows. Put for $H\in\mathcal{J}\left(  \mu\right)  $\
\[
\lambda_{E,F}(H)=(H\cap E)+F;\text{ \ }\rho_{E,F}(H)=(H+F)\cap E.
\]

The following result is obvious.

\begin{proposition}
\textit{For any }$W\in\mathcal{J}\left(  \mu\right)  $
\[
\lambda_{E,F}(W)\in\lbrack E,F];\text{ \ \ }\rho_{E,F}(W)\in\lbrack E,F];
\]

\textit{For any} $H\in\lbrack E,F]$
\[
\lambda(H)=H;\text{ \ \ }\rho(H)=H.
\]
\end{proposition}

So $\lambda_{E,F}$ and $\rho_{E,F}$ may be regarded as projections of
$\mathcal{J}\left(  \mu\right)  $ on the interval $[E,F]$. It is clear that
$\lambda^{2}=\rho\lambda=\lambda$; $\rho^{2}=\lambda\rho=\rho$.

Recall that lattices, which have the property $\lambda=\rho$, are said to be \textit{modular.}

\begin{theorem}
The lattice $\left\langle \mathcal{J}\left(  \mu\right)  ,\vee,\wedge
\right\rangle $ is modular.
\end{theorem}

\begin{proof}
From ($H\cap E)\subset^{1}\rho(H)$ and $F\subset^{1}\rho(H)$ it follows that
for all $H\in\mathcal{J}\left(  \mu\right)  $
\[
(H\cap E)+F\subset^{1}(H+F)\cap E,
\]
i.e., $\lambda(H)\subset^{1}\rho(H)$ for all $H\in\mathcal{J}\left(
\mu\right)  $.

Let $x\in(H+F)\cap E$. Recall that $[E,F]$ is an order segment, i.e.
$E\subset^{1}F$. Then $x\in H+F$ and, hence, $x=u+v$, where $v\in H$ and $u\in
F$. Its norm $\left\|  x\right\|  _{H+F}=\inf\{\left\|  v\right\|
_{H}+\left\|  u\right\|  _{F}\}$. Hence, the norm $\left\|  x\right\|  ^{1}$
of $x$, which is regarded as an element of $(H+F)\cap E$, is equal to
\[
\left\|  x\right\|  ^{1}=\max\{\inf\{\left\|  v\right\|  _{H}+\left\|
u\right\|  _{F}:v+u=x\};\text{ \ \ }\left\|  v+u\right\|  _{E}\}.
\]
The norm $\left\|  x\right\|  ^{2}$ of $x$, when $x$ is regarded as element of
$(H\cap E)+F$, is equal to
\[
\left\|  x\right\|  ^{2}=\inf\{\left\|  u\right\|  _{F},\max\{\left\|
v\right\|  _{H},\left\|  v\right\|  _{E}\}:u+v=x\}.
\]

Clearly,
\begin{align*}
\left\|  x\right\|  ^{1} &  =\max\{\inf\{\left\|  v\right\|  _{H}+\left\|
u\right\|  _{F};\text{ \ }\left\|  x\right\|  _{E}\}:v+u=x\}\\
&  \geq\max\{\inf\{\left\|  v\right\|  _{H}+\left\|  u\right\|  _{F};\text{
\ }\left\|  v\right\|  _{E}\}:v+u=x\}\\
&  \geq\inf\{\left\|  u\right\|  _{F},\max\{\left\|  v\right\|  _{H},\text{
\ }\left\|  v\right\|  _{E}\}:u+v=x\}=\left\|  x\right\|  ^{2}.
\end{align*}

Hence, $(H+F)\cap E\subset^{1}(H\cap E)+F$ and, consequently, $\lambda
_{E,F}=\rho_{E,F}$ for any interval $[E,F]$.
\end{proof}

From the property of $\left\langle \mathcal{J}\left(  \mu\right)  ,\vee
,\wedge\right\rangle $ to be modular follows \textbf{the first theorem of uniqueness.}

\begin{theorem}
Let $E$, $F\in\mathcal{J}\left(  \mu\right)  $; $E\subset^{1}F$. If there
exists such $G\in\mathcal{J}\left(  \mu\right)  $ that $G\cap E=G\cap F$ and
$G+E=G+F$ then $E=F$.
\end{theorem}

\begin{proof}
Let $G\cap E=G\cap F=X$ and $G+E=G+F=Y$.

Then $X\subset^{1}E\subset^{1}F\subset^{1}Y$ and $\left(  G\cap F\right)
+E=E\subset^{1}\left(  G+E\right)  \cap F=F$. If $E\neq F$ then $\mathcal{J}%
\left(  \mu\right)  $ is not modular.
\end{proof}

In fact, the more powerful result is true.

Recall that a lattice $A$ is said to be \textit{distributive} if for any $a$,
$b$, $c\in A$ the following equalities hold:
\[
a\vee\left(  b\wedge c\right)  =\left(  a\vee b\right)  \wedge\left(  a\vee
c\right)  ;\text{ \ }a\wedge\left(  b\vee c\right)  =\left(  a\wedge b\right)
\vee\left(  a\wedge c\right)  .
\]
It is well known that these equalities are not independent; any of them is a
consequence of other one. Both of them are equivalent to the inequality: for
any $a$, $b$, $c\in A$
\[
\left(  a\vee b\right)  \wedge c\leq a\vee\left(  b\wedge c\right)
\]
($a\leq b$ means that $a\wedge b=a$). From the last inequality it follows that
every distributive lattice is modular.

\begin{theorem}
The lattice $\left\langle \mathcal{J}\left(  \mu\right)  ,\vee,\wedge
\right\rangle $ is distributive.
\end{theorem}

\begin{proof}
It is sufficient to show that for any $E$, $F$, $G\in\mathcal{J}\left(
\mu\right)  $ the following inequality holds:
\[
\left(  E+F\right)  \cap G\subset^{1}E+\left(  F\cap G\right)  .
\]

Let $w\in\left(  E+F\right)  \cap G$. Its norm is
\begin{align*}
\left\|  w\right\|  ^{1} &  =\max\{\inf\{\left\|  u\right\|  _{E}+\left\|
v\right\|  _{F}:u+v=w\},\left\|  w\right\|  _{G}\}\\
&  =\max\{\left\|  w\right\|  _{E+F},\left\|  w\right\|  _{G}\}.
\end{align*}

Assume that $\left\|  w\right\|  _{E+F}\leq\left\|  w\right\|  _{G}$. Since
$\left\|  w\right\|  _{E+G}\leq\left\|  w\right\|  _{G}$, the norm $\left\|
w\right\|  ^{2}$ of the same element $w\in E+\left(  F\cap G\right)  $ may be
estimated as follows:
\begin{align*}
\left\|  w\right\|  ^{2} &  =\inf\{\left\|  u\right\|  _{E}+\max\{\left\|
v\right\|  _{F},\left\|  v\right\|  _{G}\}:u+v=w\}\\
&  \leq\inf\{\max\{\left\|  u\right\|  _{E}+\left\|  v\right\|  _{F},\left\|
u\right\|  _{E}+\left\|  v\right\|  _{G}\}:u+v=w\}\\
&  \leq\max\{\left\|  w\right\|  _{E+F},\left\|  w\right\|  _{E+G}%
\}\leq\left\|  w\right\|  _{G}=\left\|  w\right\|  ^{1}.
\end{align*}

If we assume that $\left\|  w\right\|  _{E+F}\geq\left\|  w\right\|  _{G}$,
then $\left\|  w\right\|  _{E+F}\geq\left\|  w\right\|  _{E+G}$ as well.
Hence,
\[
\max\{\left\|  w\right\|  _{E+F},\left\|  w\right\|  _{E+G}\}=\left\|
w\right\|  _{E+F}=\left\|  w\right\|  ^{1}.
\]
Certainly, this implies the desired inequality.
\end{proof}

As a corollary we have \textbf{the second theorem of uniqueness}.

\begin{theorem}
Let $E$, $F$, $G\in\mathcal{J}\left(  \mu\right)  $ be such that $E\cap
G=F\cap G$; $E+G=F+G$. Then either $G=E$ or $F=E$ or $G=F$.
\end{theorem}

\begin{proof}
Let $E\neq G$; $F\neq G$ and $E+G=F+G$; $E\cap G=F\cap G$. Then from
distributivity it follows that
\[
E=E\cap\left(  F+G\right)  =\left(  E\cap F\right)  +\left(  E\cap G\right)
=E\cap F.
\]

Hence $E\subset^{1}F$ and either $E=F$ or $\mathcal{J}\left(  \mu\right)  $ is
not modular (cf. theorem 4). Since every distributive lattice is modular,
$E=F$.
\end{proof}

\begin{corollary}
A pair of BIS $E$ and $F$ is uniquely determined by their sum $E+F$ and
intersection $E\cap F$.
\end{corollary}

\begin{remark}
According to the M. Stone's theorem $\left[  5\right]  $ every distributive
lattice $A$ is isomorphic (as lattice) to a some ring of sets. Moreover, as
this ring of sets it may be chosen the ring of compact open sets of the
so-called Stonian space $\frak{S}\left(  A\right)  $ of the lattice $A$ -- the
topological $T_{0}$-space, which is uniquely (up to a homeomorphism)
determined by $A$ and has the following properties:

\begin{itemize}
\item  The base of open sets of $\frak{S}\left(  A\right)  $ forms compact
open sets;

\item  Intersection of two compact open sets is compact;

\item  If $K\subset\frak{S}\left(  A\right)  $ is closed then $\cap
\{U_{i}:i\in I\}\cap K\neq\varnothing$ for an arbitrary set of compact open
sets $\{U_{i}:i\in I\}$ ($I\neq\varnothing$) of $\frak{S}\left(  A\right)  $
so that

\begin{itemize}
\item  For any $i$, $j\in I$ there is $l\in I$ such that $U_{l}\subset
U_{i}\cap U_{j}$;

\item $U_{i}\cap K\neq\varnothing$ for all $i\in I$.
\end{itemize}
\end{itemize}

If, in addition, $A$ has the maximal element, then $\frak{S}\left(  A\right)
$ is compact itself.
\end{remark}

E.g., in the aforementioned case $\left\langle \mathcal{J}\left(
\mathbb{P}\right)  ,\vee,\wedge\right\rangle $, the Stonian space
$\frak{S}\left(  \mathcal{J}\left(  \mathbb{P}\right)  \right)  $ is compact.
Other examples may be given by using Banach symmetric spaces (see the
concluding section).

\section{Closure operators and Galois connexions on $\left\langle
\mathcal{J}\left(  \mu\right)  ,\subset^{1}\right\rangle $}

Recall some algebraic definitions.

\begin{definition}
Let $\left\langle L,<\right\rangle $ be a lattice. A mapping $\pi:L\rightarrow
L$ is said to be a closure operator, if for all $a$, $b\in L$

\begin{itemize}
\item $a<b$ implies that $\pi\left(  a\right)  <\pi\left(  b\right)  $;

\item $a<\pi\left(  a\right)  $;

\item $\pi\circ\pi\left(  a\right)  =\pi\left(  a\right)  $.
\end{itemize}
\end{definition}

\begin{definition}
(Cf. $[6]$). Let $\left\langle L,<\right\rangle $ and $\left\langle L^{\prime
},<^{\prime}\right\rangle $ be lattices; $k:L\rightarrow L^{\prime}$ and
$k^{\prime}:L^{\prime}\rightarrow L$ be mappings. The pair $\left(
k,k^{\prime}\right)  $ is said to be the Galois connexion between $L$ and
$L^{\prime}$ if

\begin{itemize}
\item $a<b\Rightarrow k(a)<^{\prime}k(b)$ for $a$, $b\in L$;

\item $a^{\prime}<^{\prime}b^{\prime}\Rightarrow k^{\prime}(a^{\prime
})<k^{\prime}(b^{\prime})$ for $a^{\prime}$, $b^{\prime}\in L^{\prime}$;

\item $k^{\prime}\circ k\left(  a\right)  <a$ for $a\in L$;

\item $k\circ k^{\prime}(a^{\prime})<^{\prime}a^{\prime}$ for $a^{\prime}\in
L^{\prime}$.
\end{itemize}
\end{definition}

Below it will be needed the following simple result.

\begin{theorem}
Let $E\in\mathcal{J}_{0}\left(  \mu\right)  $. Then

\begin{enumerate}
\item $\left(  E_{0}\right)  ^{\prime\prime}\subset^{1}E^{\prime\prime}$;

\item $E^{\prime\prime}\subset^{1}\left(  \left(  E^{\prime}\right)
_{0}\right)  ^{\prime}$;

\item $\left(  \left(  \left(  E_{0}\right)  ^{\prime}\right)  _{0}\right)
^{\prime}\subset^{1}\left(  \left(  E^{\prime}\right)  _{0}\right)  ^{\prime}
$;

\item $\left(  E_{0}\right)  ^{\prime\prime}\subset^{1}\left(  \left(  \left(
E_{0}\right)  ^{\prime}\right)  _{0}\right)  ^{\prime}$.
\end{enumerate}
\end{theorem}

\begin{proof}
Since we assume that $E\in\mathcal{J}_{0}\left(  \mu\right)  $, both $E_{0}$
and $\left(  E^{\prime}\right)  _{0}$ are nontrivial. So,

1. $E_{0}\subset^{1}E\Rightarrow E^{\prime}\subset^{1}\left(  E_{0}\right)
^{\prime}\Rightarrow\left(  E_{0}\right)  ^{\prime\prime}\subset^{1}%
E^{\prime\prime}$.

2. $\left(  E^{\prime}\right)  _{0}\subset^{1}E^{\prime}\Rightarrow
E^{\prime\prime}\subset^{1}\left(  \left(  E^{\prime}\right)  _{0}\right)
^{\prime}$.

3. $E_{0}\subset^{1}E\Rightarrow E^{\prime}\subset^{1}\left(  E_{0}\right)
^{\prime}\Rightarrow\left(  E^{\prime}\right)  _{0}\subset^{1}\left(  \left(
E_{0}\right)  ^{\prime}\right)  _{0}\Rightarrow\left(  \left(  \left(
E_{0}\right)  ^{\prime}\right)  _{0}\right)  ^{\prime}\subset^{1}\left(
\left(  E^{\prime}\right)  _{0}\right)  ^{\prime}$.

4. $\left(  \left(  E_{0}\right)  ^{\prime}\right)  _{0}\subset^{1}\left(
E_{0}\right)  ^{\prime}\Rightarrow\left(  E_{0}\right)  ^{\prime\prime}%
\subset^{1}\left(  \left(  \left(  E_{0}\right)  ^{\prime}\right)
_{0}\right)  ^{\prime}$.
\end{proof}

Consider the lattice $\mathcal{J}_{0}^{\ast}\left(  \mu\right)  $, equipped
with the inverse order $\vartriangleleft$: $E\vartriangleleft F$ is the same
as $F\subset^{1}E$ for $E$, $F\in\mathcal{J}_{0}\left(  \mu\right)  $. Put
$\mathcal{J}_{0}:=\left\langle \mathcal{J}_{0}\left(  \mu\right)  ,\subset
^{1}\right\rangle $; $\mathcal{J}_{0}^{\ast}:=\left\langle \mathcal{J}%
_{0}\left(  \mu\right)  ,\vartriangleleft\right\rangle $; $\mathcal{J}%
:=\left\langle \mathcal{J}\left(  \mu\right)  ,\subset^{1}\right\rangle $

\begin{theorem}
The mapping $\left(  _{0}\right)  :\mathcal{J}_{0}^{\ast}\left(  \mu\right)
\rightarrow\mathcal{J}_{0}^{\ast}\left(  \mu\right)  $ is a closure operator
on $\mathcal{J}_{0}^{\ast}$.

The mapping $\left(  ^{\prime\prime}\right)  :\mathcal{J}_{0}\left(
\mu\right)  \rightarrow\mathcal{J}_{0}\left(  \mu\right)  $ is a closure
operator on $\mathcal{J}$.
\end{theorem}

\begin{proof}
The proof is an obvious consequence of definitions and the summary.
\end{proof}

It may be defined the most important sublattices of the $\mathcal{J}\left(
\mu\right)  $

\begin{definition}
The lattice $\mathcal{J}_{00}\left(  \mu\right)  $ consists of all BIS $E$
with the absolute continuous norm (i.e., such that $E=E_{0}$).

The lattice $\mathcal{J}^{\prime}\left(  \mu\right)  $ contains all BIS $E$ of
kind $E=F^{\prime}$ for some BIS $F$.
\end{definition}

\begin{corollary}
Lattices $\mathcal{J}_{00}\left(  \mu\right)  $ and $\mathcal{J}^{\prime
}\left(  \mu\right)  $ both are distributive (and, hence, modular). If the
measure $\mu$ is a probability measure $\mathbb{P}$ and all BIS $E$ from
$\mathcal{J}\left(  \mathbb{P}\right)  $ satisfy the norming condition (i.e.
if $\mathcal{J}\left(  \mathbb{P}\right)  $ is Dedekind complete) then
sublattices $\mathcal{J}_{00}\left(  \mathbb{P}\right)  $ and $\mathcal{J}%
^{\prime}\left(  \mathbb{P}\right)  $ both are Dedekind complete too.
\end{corollary}

\begin{proof}
As it is well known, the set of fixed points of a closure operator that acting
on the distributive Dedekind complete lattice $A$ is a sublattice of $A$ that
holds these properties. Clearly, the set of fixed points of $\left(
_{0}\right)  $ is exactly $\mathcal{J}_{00}\left(  \mu\right)  $. To show that
$\mathcal{J}^{\prime}\left(  \mu\right)  $ coincides with the set of fixed
points of the closure operator $\left(  ^{\prime\prime}\right)  $ it is enough
to notice that $E=F^{\prime}$ for some BIS $F$ if and only if $E=E^{\prime
\prime}$. Certainly, if $E=E^{\prime\prime}$ then $E=F^{\prime}$ for
$F=E^{\prime}$. Conversely, if $E=F^{\prime}$ then $E^{\prime\prime}%
=F^{\prime\prime\prime}=\left(  F^{\prime}\right)  ^{\prime\prime}=F^{\prime}$.
\end{proof}

Correspondences $E\rightsquigarrow E_{0}$ and $E\rightsquigarrow E^{\prime}$
may be regarded as mappings. It will be written $\left(  _{0}\right)  :$
$E\rightarrow E_{0}$; $\left(  ^{\prime}\right)  :E\rightarrow E^{\prime}$.

Using the mappings $\left(  _{0}\right)  $ and $\left(  ^{\prime}\right)  $ it
may be constructed a pair of mappings, say $k$ and $k^{\prime}$ that are in
the Galois connexion.

Let $k:\mathcal{J}_{0}\left(  \mu\right)  \rightarrow\mathcal{J}_{0}^{\ast
}\left(  \mu\right)  $ and $k^{\prime}:\mathcal{J}_{0}^{\ast}\left(
\mu\right)  \rightarrow\mathcal{J}_{0}\left(  \mu\right)  $ are given by
\[
kE=\left(  E_{0}\right)  ^{\prime};\text{ \ }k^{\prime}E=\left(  E^{\prime
}\right)  _{0}.
\]

\begin{theorem}
The pair $\left(  k,k^{\prime}\right)  $ is the Galois connexion between
$\mathcal{J}_{0}^{\ast}$ and $\mathcal{J}_{0}$.
\end{theorem}

\begin{proof}
Let us check up properties from the definition 5.

\begin{itemize}
\item $E\subset^{1}F\Rightarrow E_{0}\subset^{1}F_{0}\Rightarrow\left(
F_{0}\right)  ^{\prime}\subset^{1}\left(  E_{0}\right)  ^{\prime}%
\Rightarrow\left(  E_{0}\right)  ^{\prime}\vartriangleleft\left(
F_{0}\right)  ^{\prime} $.

\item $E\vartriangleleft F\Rightarrow F\subset^{1}E\Rightarrow\left(
E^{\prime}\right)  _{0}\subset^{1}\left(  F^{\prime}\right)  _{0}$.

\item $((\left(  F^{\prime}\right)  _{0})_{0})^{\prime}=(\left(  F^{\prime
}\right)  _{0})^{\prime}$. By the theorem 4, $F^{\prime\prime}\subset
^{1}(\left(  F^{\prime}\right)  _{0})^{\prime}$ and, hence,
\[
F\subset^{1}F^{\prime\prime}\subset^{1}(\left(  F^{\prime}\right)
_{0})^{\prime}.
\]

\item $(\left(  \left(  E_{0}\right)  ^{\prime}\right)  ^{\prime})_{0}=\left(
\left(  E_{0}\right)  ^{\prime\prime}\right)  _{0}\subset^{1}E$. Hence,
$E\vartriangleleft\left(  \left(  E_{0}\right)  ^{\prime\prime}\right)  _{0}$.
\end{itemize}
\end{proof}

\begin{corollary}
Compositions $k\circ k^{\prime}$and $k^{\prime}\circ k$\ are closure operators
on $\mathcal{J}_{0}^{\ast}$ and $\mathcal{J}_{0}$ respectively.
\end{corollary}

\begin{proof}
This is an obvious consequence of definitions.\ 
\end{proof}

\begin{remark}
So, we obtain some more closure operators on $\mathcal{J}_{0}\left(
\mu\right)  $. Notice that $k\circ k^{\prime}:E\rightarrow(\left(  E^{\prime
}\right)  _{0})^{\prime}$ and $k^{\prime}\circ k:E\rightarrow\left(  \left(
E_{0}\right)  ^{\prime\prime}\right)  _{0}$. The second mapping is coincide
with the usual $\left(  _{0}\right)  :E\rightarrow E_{0}$. However the first
one pick out from $\mathcal{J}^{\prime}\left(  \mu\right)  $ those BIS that
are dual to BIS having the absolutely continuous norm.
\end{remark}

\section{Symmetric Banach spaces}

Results from previous sections are of the most interest when a special class
of BIS - the class of \textit{symmetric Banach spaces} is considered.

Recall the definition.

\begin{definition}
A Banach ideal space $E$ of (classes of) measurable real functions, which are
defined on the admissible measure space $\left(  \Omega,\Sigma,\mu\right)  $
is said to be \textit{symmetric }if for any functions $x=x\left(  t\right)  $
and $y=y\left(  t\right)  $ of $E$ the following condition holds:

\begin{itemize}
\item  If $x\in E$ and functions $\left|  y(t)\right|  $ and $\left|
x(t)\right|  $ are equimeasurable then $y\in E$ and $\left\|  y\right\|
_{E}=\left\|  x\right\|  _{E}$.
\end{itemize}
\end{definition}

Let $\mathcal{S}\left(  \Omega,\Sigma,\mu\right)  $ be a class of all
symmetric Banach spaces (they in the future will be referred to as symmetric
spaces). This class contains Lebesgue-Riesz spaces $L_{p}\left(  \mu\right)
$; Orlicz spaces $L_{M}\left(  \mu\right)  $ and so on. Usually properties 1
and 2 are supplemented with the following \textit{norming condition}:
\[
\left\|  \chi_{e}\right\|  _{E}=1\text{ \ for any set }e\in\Sigma\text{ \ of
the measure }\mu\left(  e\right)  =1.
\]

This condition implies that $\mathcal{S}\left(  \Omega,\Sigma,\mu\right)  $ is
a Dedekind complete distributive lattice because of the known theorem of
inclusion:
\[
L_{1}\left(  \mu\right)  \cap L_{\infty}\left(  \mu\right)  \subset
^{1}E\left(  \mu\right)  \subset^{1}L_{1}\left(  \mu\right)  +L_{\infty
}\left(  \mu\right)  .
\]
For a probability measure $\mathbb{P}$ these inclusions looks like
\[
L_{\infty}\left(  \mathbb{P}\right)  \subset^{1}E\left(  \mathbb{P}\right)
\subset^{1}L_{1}\left(  \mathbb{P}\right)  .
\]

For a purely atomic measure (with mass of every point is equal to 1) we obtain
so called \textit{symmetric discrete (or sequence) spaces }defined on an
arbitrary set. The only characteristic that distinguishes corresponding
classes of discrete spaces is the cardinality of $\Omega$. For
$\operatorname{card}\Omega=\varkappa$ \ the class $\mathcal{S}\left(
\Omega,\Sigma,\mu\right)  $ will be denoted by $\mathcal{S}\left(
\varkappa\right)  $. Inclusions in this case looks like:
\[
l_{1}\left(  \varkappa\right)  \subset^{1}E\left(  \varkappa\right)
\subset^{1}l_{\infty}\left(  \varkappa\right)  .
\]

Notice that in the case of symmetric spaces the finiteness of $\mu$ implies
that it is non-atomic.

As it follows from the preceding consideration, all lattices of symmetric
Banach spaces may be participate into three parts:

\begin{itemize}
\item $\mathcal{S}^{\left(  1\right)  }$ - lattices of symmetric spaces,
defined on a probability (non atomic) space;

\item $\mathcal{S}^{\left(  \infty\right)  }$ - lattices of symmetric spaces,
defined on a non atomic space of infinite measure;

\item $\mathcal{S}^{\left(  \mathbf{D}\right)  }$ - lattices of symmetric
sequence spaces:
\[
\mathcal{S}=\mathcal{S}^{\left(  1\right)  }\cup\mathcal{S}^{\left(
\infty\right)  }\cup\mathcal{S}^{\left(  \mathbf{D}\right)  }.
\]
\end{itemize}

Our nearest aim is to show that all lattices $\left\langle \mathcal{S}\left(
\mu\right)  ,\subset^{1}\right\rangle $ from a given class $\mathcal{S}%
^{\left(  \mathbf{?}\right)  }$, where $\mathbf{?}\in\{1,\infty,\mathbf{D}\}$
are pairwice lattice isomorphic, i.e. that there are at most three different
lattices amongst all of kind $\mathcal{S}\left(  \mu\right)  .$

\begin{theorem}
Lattices $\mathcal{S}\left(  \mu\right)  $ and $\mathcal{S}\left(  \nu\right)
$ that belong to the same class $\mathcal{S}^{\left(  \mathbf{?}\right)  }$,
where $\mathbf{?}\in\{1,\infty,\mathbf{D}\}$ are lattice-isomorphic.
\end{theorem}

\begin{proof}
The one-to-one correspondence between members of these lattices may be shown
by using the operation of replanting of symmetric spaces from one measure to
another. Such operation was suggested by A.A. Mekler [7].

Namely, let $E\left(  \mu\right)  \in$ $\mathcal{S}\left(  \mu\right)  $. To
$x\left(  t\right)  \in E\left(  \mu\right)  $ corresponds its
\textit{distribution function}
\[
n_{x}\left(  s\right)  =\operatorname{mes}\left(  \{t\in\lbrack0,1]:x\left(
t\right)  >s\}\right)
\]
and its \textit{non-increasing rearrangement}
\[
x^{\ast}\left(  t\right)  =\inf\{s\in\lbrack0,\infty):n_{\left|  x\right|
}\left(  s\right)  <t\}.
\]

Obviously, $x^{\ast}\left(  t\right)  $ is defined either on $\left[
0,1\right]  $ or on $[0,\infty)$ (both with the Lebesgue measure) or at
$\mathbb{N}$ (with mass 1 in every point) and, hence is an element of the
corresponding vector lattice $L_{0}\left[  0,1\right]  $ (resp., $L_{0}\left[
0,\infty\right]  $) or $L_{0}\left[  \mathbb{N}\right]  $. Notice that usually
$L_{0}\left[  \mathbb{N}\right]  $ is denoted by $s$).

Let, for distinctness, $x^{\ast}\left(  t\right)  \in L_{0}\left[  0,1\right]
$.

It is clear that the set
\[
\widetilde{E}\left[  0,1\right]  =\{x^{\ast}\left(  \sigma t\right)  :x\left(
t\right)  \in E\left(  \mu\right)  ;\text{ }\sigma\in\operatorname{Aut}\}
\]
(where $\operatorname{Aut}$ denotes the set of all preserving measure
automorphisms of $\left[  0,1\right]  $) is an ideal in $L_{0}\left[
0,1\right]  $, which, being equipped with the norm
\[
\left\|  x^{\ast}\left(  t\right)  \right\|  _{\widetilde{E}}:=\left\|
x\right\|  _{E\left(  \mu\right)  }\text{; \ }\left\|  y\left(  t\right)
\right\|  _{\widetilde{E}}:=\left\|  y^{\ast}\left(  t\right)  \right\|
_{\widetilde{E}},
\]
\ becomes a Banach symmetric function space on $\left[  0,1\right]  $.\ 

Moreover, this space is uniquely determined by $E\left(  \mu\right)  $.

Notice that for any probability measure $\nu$ the space $\widetilde{E}\left[
0,1\right]  $ in the same way uniquely defines the corresponding space
$\widetilde{E}\left(  \nu\right)  $.

It will be said that the space $\widetilde{E}\left(  \nu\right)  $ is obtained
from $E\left(  \mu\right)  $ by the\textit{\ replanting Mekler's
procedure:}$\ $%
\[
\operatorname{Mekl}_{\mu,\nu}:\mathcal{S}\left(  \mu\right)  \rightarrow
\mathcal{S}\left(  \nu\right)  ;
\]%
\[
E\left(  \mu\right)  \rightarrowtail\widetilde{E}\left(  \nu\right)
\]

So, every symmetric space $E=E\left(  \mu\right)  $ generates a tower
\[
\left\lfloor E\right\rfloor =\{\widetilde{E}\left(  \Omega,\Sigma,\nu\right)
=\operatorname{Mekl}_{\mu,\nu}(E(\mu))\}
\]
where $\left(  \Omega,\Sigma,\nu\right)  $ runs all probability spaces.

Obviously, the procedure $\operatorname{Mekl}_{\mu,\nu}$ holds the relation
$\subset^{1}$and lattice operations:%

\[
E\subset^{1}F\Rightarrow\operatorname{Mekl}_{\mu,\nu}(E)\subset^{1}%
\operatorname{Mekl}_{\mu,\nu}(F);
\]%
\begin{align*}
\operatorname{Mekl}_{\mu,\nu}(E\cap F) &  =\operatorname{Mekl}_{\mu,\nu
}(E)\cap\operatorname{Mekl}_{\mu,\nu}(F);\\
\operatorname{Mekl}_{\mu,\nu}(E+F) &  =\operatorname{Mekl}_{\mu,\nu
}(E)+\operatorname{Mekl}_{\mu,\nu}(F)\text{;}%
\end{align*}%
\[
\operatorname{Mekl}_{\mu,\nu}(E_{0})=\left(  \operatorname{Mekl}_{\nu
}(E)\right)  _{0};\text{ \ \ }\operatorname{Mekl}_{\mu,\nu}(E^{\prime
})=\left(  \operatorname{Mekl}_{\mu,\nu}(E)\right)  ^{\prime}%
\]

and, hence, generates the lattice isomorphism between lattices $\mathcal{S}%
\left(  \mu\right)  $ and $\mathcal{S}\left(  \nu\right)  $.

Similarly for lattices that are belong to $\mathcal{S}^{\left(  \mathbf{?}%
\right)  }$, where either $\mathbf{?}=\infty$, or $\mathbf{?}=\mathbf{D}$.
\end{proof}

\begin{remark}
It is not clear: whether are lattice-isomorphic lattices of different types
(e.g., from $\mathcal{S}^{\left(  1\right)  }$and $\mathcal{S}^{\left(
\infty\right)  }$). It may be shown that lattices $\mathcal{S}^{\left(
1\right)  }$ and $\mathcal{S}^{\left(  \mathbf{D}\right)  }$ are isomorphic to
quotient lattices of $\mathcal{S}^{\left(  \infty\right)  }$
\end{remark}

\begin{summary}
There exists (up to lattice-isomorphism) at most three lattices $\mathcal{S}%
=\mathcal{S}^{\left(  \mathbf{?}\right)  }$ ($?\in\{1,\infty,\mathbf{D}\}$) of
symmetric spaces.

All of them are Dedekind complete, distributive (and, hence, modular).

Operations $kE=\left(  E_{0}\right)  ^{\prime}\ $and \ $k^{\prime}E=\left(
E^{\prime}\right)  _{0}$ interrelate the ''fundamental part'' $\mathcal{S}%
_{0}$ (see the definition 3) of each of these lattices with itself in the
reverse order (i.e., with $\mathcal{S}_{0}^{\ast}$ by the given notation):

The pair $\left(  k.k^{\prime}\right)  $ is the Galois connexion between
$\mathcal{S}_{0}^{\ast}$ and $\mathcal{S}_{0}$
\end{summary}

\begin{remark}
For the class $\mathcal{S}$ of all symmetric spaces, classes $\mathcal{S}%
_{00}$ and $\mathcal{S}^{\prime}$ (see the definition 6) are just those, whose
elements are usually called rearrangement invariant Banach spaces
\end{remark}

\section{References}

\begin{enumerate}
\item  Buchvalov A.V., Veksler A.I., Geiler V.A. \textit{Normed lattices} (in
Russian), Math. Analysis, (Itogi Nauki i Techniki, VINITI) \textbf{18 }(1980) 125-184

\item  Buchvalov A.V., Veksler A.I., Losanovskii G.Ya. \textit{Banach lattices
- some Banach aspects of the theory }(in Russian), Uspechi Mat. Nauk,
\textbf{34:2} (1979) 179-484

\item  Seever G.L. \textit{A peculiar Banach function spaces}, Proc. AMS
\textbf{16} (1965) 662-664

\item  Aronszajn N., Gagliardo E. \textit{Interpolation spaces and
interpolation methods}, Ann. Math. Pura App. Ser. 4, \textbf{68} (1965) 51-117

\item  Stone M.H. \textit{Topological representations of distributive lattices
and Browerian logics}, \v{C}asopis pro P\u{e}st. Mat. Fyz. \textbf{6}7 (1937) 1-25

\item  Ore O. \textit{Galois connexions}, Trans. AMS \textbf{55} (1944) 493-513

\item  Mekler A.A. \textit{On one equivalence relation in the class of
symmetric ideals of the space }$L_{1}$ (in Russian), Proc. \textbf{VI}
All-Union Seminar on superfluity problems and informational systems, Part
\textbf{1}, Leningrad, (1974) 104-106
\end{enumerate}
\end{document}